\documentclass{llncs}
\usepackage{makeidx}  % allows for indexgeneration
\usepackage{graphicx}
\begin{document}
\pagestyle{headings}     % switches on printing of running heads
\mainmatter              % start of the contributions
\title{Forbidden substrings, Kolmogorov complexity
and almost periodic sequences}
\titlerunning{Forbidden substrings, Kolmogorov complexity
and almost periodic sequences}   %    abbreviated title (for running head)
\author{A.~Yu.~Rumyantsev, M.~A.~Ushakov}
\authorrunning{A.~Yu.~Rumyantsev, M.~A.~Ushakov}   % abbreviated author list (for running head)
%
%%%% modified list of authors for the TOC (add the affiliations)
\tocauthor{A.~Yu.~Rumyantsev, M.~A.~Ushakov (Moscow State University)}
\institute{Moscow State University, Russia, Mathematics Department,
Logic and algorithms theory division}

\maketitle              % typeset the title of the contribution

\begin{abstract}
Assume that for some $\alpha<1$ and for all nutural~$n$ a set $F_n$ of at
most $2^{\alpha n}$ ``forbidden'' binary strings of length~$n$ is fixed.
Then there exists an infinite binary sequence $\omega$ that does not
have (long) forbidden substrings.

We prove this combinatorial statement by translating it into a
statement about Kolmogorov complexity and compare this proof
with a combinatorial one based on Laslo Lovasz local lemma.

Then we construct an almost periodic sequence with the same
property (thus combines the results from~\cite{dls} and~\cite{msu}).

Both the combinatorial proof and Kolmogorov complexity argument
can be generalized to the multidimensional case.
\end{abstract}

\section{Forbidden strings}

Fix some positive constant $\alpha<1$. Assume that for each
natural $n$ a set $F_n$ of binary strings of length $n$ is
fixed. Assume that $F_n$ consists of at most $2^{\alpha n}$
strings.

We look for an infinite binary sequence $\omega$ that does not
contain forbidden substrings.

\begin{proposition}\label{forbidden-strings}
š š š š %
There exists an infinite binary sequence $\omega$ and a constant
$N$ such that for any $n>N$ the sequence $\omega$ does not have
a substring $x$ of length~$n$ that belongs to $F_n$.
š š š š %
\end{proposition}

One may consider strings in $F_n$ as ``forbidden'' strings of
length~$n$; proposition then says that there exists an infinite
sequence without (sufficiently long) forbidden substrings.

For example, we can forbid strings having low Kolmogorov
complexity. Let $F_n$ be the set of all strings of
length~$n$ whose complexity is less than $\alpha n$. Then
$\#F_n$ does not exceed $2^{\alpha n}$ (there are at most
$2^{\alpha n}$ programs of size less than $\alpha n$).

Therefore Proposition~\ref{forbidden-strings} implies the
following statement that was used in~\cite{dls}:

\begin{proposition}\label{levin-lemma}
š š š š %
For any $\alpha<1$ there exists a number $N$ and an infinite
binary sequence $\omega$ such that any its substring $x$ of
length greater than $N$ has high complexity:
š š š š $$
{\mathrm{K}}(x) \ge \alpha |x|.
š š š š $$
\end{proposition}

Here ${\mathrm{K}}(x)$ stands for Kolmogorov complexity of $x$ (the length
of the shortest program producing $x$,
the definition is given in~\cite{livitan});
it does not matter which version of Kolmogorov compexity
(prefix, plain, etc.) we consider since the logarithmic difference
between them can be compensated by a small change in~$\alpha$.
The notation $|x|$ means the length of string $x$.

Our observation is that the reverse implication is true, i.e.,
Proposition~\ref{levin-lemma} implies Proposition~\ref{forbidden-strings}.
It is easy to see if we consider a stronger version of
Proposition~\ref{levin-lemma} when ${\mathrm{K}}$ is replaced
by a relativized version ${\mathrm{K}}_A$ where $A$
is an arbitrary oracle (an external procedure that can be
called). Indeed, consider the set of all forbidden strings as an
oracle. Then the relativized complexity of any string in $F_n$
does not exceed $\alpha n + O(1)$ since its ordinal number in
the length-sorted list of all forbidden strings is at most
$\sum_{k\le n}2^{\alpha k}=O(2^{\alpha n})$. The constant $O(1)$
can be absorbed by a small change in~$\alpha$, and we get the
statement of Proposition~\ref{forbidden-strings}.

More interestingly, we can avoid relativization and derive
Proposition~\ref{forbidden-strings} from (non-relativized)
Propostion~\ref{levin-lemma}. It can be done as follows.

First note that we may assume (without loss of generality) that
$\alpha$ is rational. Assume that for some set $F$ of forbidden
strings the statement of Proposition~\ref{forbidden-strings} is
false. Then for each $c\in\bbbn$ there exists a set $F^c$
with the following properties:

(a)~$F^c$ consists of strings of length greater than~$c$;

(b)~$F^c$ contains at most $2^{\alpha k}$ strings of length $k$
for any $k$;

(c)~any infinite binary string has at least one substring that
belongs to~$F^c$.

(Indeed, let $F^c$ be the set of all forbidden strings
that have length more than~$c$.)

The statement (c) can be reformulated as follows: the family of
open sets $S_x$ for all $x\in F^c$ covers the set $\Omega$ of
all binary sequences, where $S_x$ is a set of all sequences that
have substring $x$. The standard compactness argument implies
that $F^c$ can be replaced by its finite subset, so we assume
without loss of generality that $F^c$ is finite.

The properties (a), (b) and~(c) are enumerable (for finite
$F^c$): each $S_x$ is an enumerable union of intervals, so if
the sets $S_x$ for $x\in F^c$ cover $\Omega$, this can be
discovered at a finite step. (In fact, they are decidable,
but this does not matter.) So the first set $F_c$ encountered in the
enumeration (for a given~$c$) is a computable function of~$c$.

Now we can construct a decidable set of forbidden strings that
does not satisfy the statement of
Proposition~\ref{forbidden-strings}. Indeed, construct a
sequence $c_1< c_2< c_3<\ldots$ where $c_{i+1}$ is greater than
the length of all strings in $F^{c_i}$ and take the union
of all $F^{c_i}$. We obtain the decidable set~$\hat F$ such that
$\hat F$ contains at most $2^{\alpha k}$ strings of length $k$
for any $k$, and any infinite binary string has (for any $i$) at
least one substring of length greater that $c_i$ that belongs
to~$\hat F$. For this decidable set we need no special oracle,
q.e.d.

The proof of Proposition~\ref{levin-lemma} given
in~\cite{dls} uses prefix complexity. See below
Section~\ref{almost-periodic} where we prove the stronger
version of this Proposition needed for our purposes.

\section{Combinatorial proof}

The statement of Proposition~\ref{forbidden-strings} has nothing
to do with Kolmogorov complexity. So it would be natural to look
for a combinatorial proof.

The simplest idea is to use the random bits as the elements of the
sequence. Then the probability of running into a forbidden string
in a given $k$ positions
š š š š $$
\omega_{n}\omega_{n+1}\ldots\omega_{n+k-1}
š š š š $$
is bounded by $2^{-(1-\alpha)k}$, i.e., exponentially decreases
when $k\to\infty$. However, the number of positions where a
forbidden string of a given length can appear is infinite, and
the sum of probablities is infinite too. And, indeed, a truly
random sequence contains any string as its substring, so we need
to use something else.

Note that two non-overlapping fragments of a random sequence are
independent. So the dependence can be localized and we can apply
the following well-known statement:

\begin{proposition}[Laslo Lovasz local lemma]
š š š š %
Let $G$ be a graph with vertex set $V=\{v_1,\ldots,v_n\}$ and
edge set $E$. Let $A_i$ be some event associated with vertex
$v_i$. Assume that for each~$i$ the event $A_i$ is independent
with the random variable ``outcomes of all $A_j$ such that $v_j$
is not connected to $v_i$ by an edge''. Let $p_i\in(0,1)$ be a
number associated with $A_i$ in such a way that
š š š š $$
\Pr [A_i] \le p_i \prod_{v_j\sim v_i} (1-p_j)
š š š š $$
where the product is taken over all neighbour vertices $v_j$
\textup(connected to $v_i$ by an edge\textup). Then
š š š š $$
\Pr[\mbox{neither of $A_i$ happens}] \ge \prod_{i=1}^n (1-p_i)
š š š š $$
and, therefore, this event is non-empty.
\end{proposition}

The proof of this Lemma could be found, e.g.,
in~\cite{randomized}, p.~115.

To apply this Lemma to our case consider a finite random string
of some fixed length~$N$ where all bits are independent and
unbiased (both outcomes have probability~$1/2$). Consider a
graph whose vertices are intervals of indices (i.e., places
where a substring is located) of length at least $L$ (some
constant to be chosen later). Two intervals are connected by an
edge if they are not disjoint (share some bit). For each
interval $v$ consider the event $A_v$: ``substring of the random
string located at~$v$ is forbidden''. This event is independent
with all events that deal with bits outside $v$, so the
independence condition is fulfilled.

Let $p_v=2^{-\delta|v|}$ for all $v$ and some $\delta$ (to be
chosen later). To apply the lemma, we need to prove that
š š š š $$
\Pr[A_v] \le p_v\prod_{\mbox{\scriptsize $v$ and $w$ are}\atop\mbox{\scriptsize not disjoint}}(1-p_w).
š š š š $$

Let $l\ge L$ be the length of the string $v$ and let
š š š š $$
R=\prod_{\mbox{\scriptsize $v$ and $w$ are}\atop\mbox{\scriptsize not disjoint}}(1-p_w).
š š š š $$
Then
š š š š $$
R\ge \prod_{k=L}^{N} (1-2^{-\delta k})^{l+k}
š š š š $$
(strings $w$ have length between $L$ and $N$ and there are at
most $l+k$ strings of length $k$ that share bits with~$v$), and
š š š š $$
R\ge \left[\prod_{k\ge L} (1-2^{-\delta k})\right]^l
š š š\prod_{k\ge L} (1-2^{-\delta k})^k
š š š š $$
(we split the product in two parts and replace finite products
by infinite ones). The product $\prod (1-\varepsilon_i)$
converges if and only if the series $\sum\varepsilon_i$ converges.
The corresponding series
š š š š $$
\sum_{k\ge L} 2^{-\delta k} \mbox{ and }
\sum_{k\ge L} k\cdot 2^{-\delta k}
š š š š $$
do converge. Therefore both products converge
and for a large $L$ both products are close to~$1$:
š š š š $$
R\ge C_1^l C_2 \ge D^l
š š š š $$
where $C_1$, $C_2$ and $D$ are some constants that could be made
close to~$1$ by choosing a large enough $L$ (not
depending on~$l$). Then
š š š š $$
p_v R \ge 2^{-\delta l} D^l \ge 2^{-\delta l} 2^{-\gamma l}=
2^{-(\delta+\gamma)l},
š š š š $$
where $\gamma=-\log D$ could be arbitrarily small for some~$L$.
We choose $\delta$ and $L$ in such a way that
$\delta<(1-\alpha)/2$ and $\gamma<(1-\alpha)/2$. Then
š š š š $$
p_v R \ge 2^{-(1-\alpha)l}\ge \Pr[A_v]
š š š š $$
(forbidden strings form a $2^{-(1-\alpha)l}$-fraction of all
strings having length~$l$) and conditions of Lovasz lemma are
fulfilled.

So we see that for some large $L$ and for all sufficiently
large~$N$ there exists a string of length~$N$ that does not
contain forbidden strings of length $L$ or more. Standard
compactness argument shows that there exists an infinite
binary string with the same property.

This finishes the combinatorial proof of
Proposition~\ref{forbidden-strings}.

Note that this combinatorial proof hardly can be considered as a
mere translation of Kolmogorov complexity argument. Another
reason to consider it as a different proof is that it has a
straightforward generalization for several dimensions. (The
Kolmogorov complexity argument has this too, as we see in
Section~\ref{multidimensional}, but requires significant
changes.)

A $d$-dimensional sequence is a function
$\omega\colon\bbbz^d\to\{0,1\}$. Instead of substrings we
consider $d$-dimensional ``subcubes'' in the sequence, i.e.,
restrictions of $\omega$ to some cube in $\bbbz^d$. For any
$n$ there are $2^{n^d}$ different cubes with side $n$. Assume that
for every $n>1$ a set $F_n$ of not more than $2^{\alpha n^d}$
``forbidden cubes'' is fixed.

\begin{proposition}\label{forbidden-strings-multidimensional}
š š š š %
There exists a number $L$ and $d$-dimensional sequence that does
not contain forbidden subcube with side greater than $L$.
š š š š %
\end{proposition}

The proof repeats the combinatorial proof of
Proposition~\ref{forbidden-strings} with the following changes.
The bound for~$R$ now is
š š š š $$
R\ge \prod_{k=L}^{N} (1-2^{-\delta k^d})^{(l+k)^d},
š š š š $$
since there are at most $(l+k)^d$ cubes with side~$k$ intersecting a
given cube with side~$l$. Then we represent $(l+k)^d$ as a sum
of $d+1$ monomials and get a representation of this bound as a
product of infinite products, each for one monomial. Every
product has the following form (for some~$i$ in $0\ldots d$ and
for some $c_i$ that depends on $d$ and $i$, but not $k$ and $l$):
š š š š $$
\prod_{k\ge L} (1-2^{-\delta k^d})^{c_i l^i k^j}
= \left[\prod_{k\ge L} (1-2^{-\delta k^d})^{k^j}\right]^{c_i l^i}.
š š š š $$
The corresponding series obviously converge (due to the same reasons as
before), and again we can make expression $[\ldots]$ as close to $1$ as needed
by choosing $L$ (and again the choice of $L$ does not depend on
$l$). Then the estimate for $R$ takes the form:
š š š š $$
R \ge \prod_{i=0}^{d} D_i^{c_i l^i} \ge \prod_{i=1}^{d} D_i^{Cl^d}
\ge \left[\prod_{i=1}^{d} D_i^C\right]^{l^d}
\ge D^{l^d},
š š š š $$
where $c_i$, $D_i$, $C$ and~$D$ are some constants, and $C$ and~$D$
could be made as close to $1$ as needed.

Then the proof goes exactly as before.

\section{Construction of almost periodic sequences}
\label{almost-periodic}

A sequence is called \emph{almost periodic} if each of its
substrings has infinitely many occurences at limited distances,
i.e., for any substring~$x$ there exists a number~$k$ such that
any substring~$y$ of $\omega$ of length~$k$ contains $x$.

The following result is proven in~\cite{msu} (in the paper almost
periodic sequences were called strongly almost periodic sequences):

\begin{proposition}\label{prefixes}
š š š š %
Let $\alpha<1$ be a constant. There exists an almost periodic
sequence~$\omega$ such that any sufficiently long prefix~$x$
of~$\omega$ has large complexity: ${\mathrm{K}}(x)\ge\alpha|x|$.
š š š š %
\end{proposition}

Comparing this statement with Proposition~\ref{levin-lemma}, we
see that there is an additional requirement for the sequence to
be almost periodic; on the other hand
high complexity is guaranteed only for
prefixes (and not for all substrings).

Now we combine these two results:

\begin{proposition}\label{main}
š š š š %
Let $\alpha<1$ be a constant. There exists an almost periodic
sequence~$\omega$ such that any sufficiently long substring~$x$
of~$\omega$ has large complexity: ${\mathrm{K}}(x)\ge\alpha|x|$.
š š š š %
\end{proposition}

The paper~\cite{msu} provides a universal construction for almost
periodic sequences. Now we suggest another, less general
construction that is more suitable for our purposes.

Namely, we define some equivalence relation on the set of
indices ($\bbbn$). Then we construct a sequence
š š š š $$
\omega=\omega_0\omega_1\omega_2\ldots
š š š š $$
with the following property: $i\equiv j \Rightarrow
\omega_i=\omega_j$. In other words, all the places that belong
to one equivalence class carry the same bit. This property
guarantees that $\omega$ is almost periodic if the equivalence
relation is chosen in a proper way.

Let $n_0, n_1, n_2,\ldots$ be an increasing sequence of natural
numbers such that $n_{i+1}$ is a multiple of $n_i$ for each $i$.
The prefix of length $n_0$, i.e., the interval $[0,n_0)$, is
repeated with period $n_1$. This means that for any $i$ such
that $0\le i < n_0$ the numbers
š š š š $$
i, i+n_1, i+2n_1, i+3n_1,\ldots
š š š š $$
belong to the same equivalence class. In the similar way the
interval $[0,n_1)$ is repeated with period $n_2$: for any $i$
such that $0\le i < n_2$ the numbers
š š š š $$
i, i+n_2, i+2n_2, i+3n_2,\ldots
š š š š $$
are equivalent. (Note that $n_2$ is a multiple of $n_1$,
therefore the equivalence classes constructed at the first step
are not changed.) And so on: for any $i\in [0,n_s)$ and for any
$k$ the numbers $i$ and $i+kn_{s+1}$ are equivalent.

\begin{figure}[h]
\begin{center}
\includegraphics[scale=0.6]{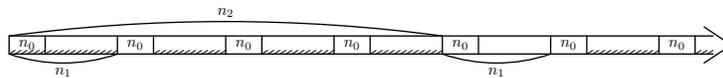}
\end{center}
\caption{Primary (shaded) and secondary bits in a sequence}\label{bits}
\end{figure}

The following statement is almost evident:

\begin{proposition}\label{construction}
š š š š %
If a sequence $\omega$ respects this equivalence relation,
i.e., the equivalent positions have equal bits, then the
sequence in almost periodic.
š š š š %
\end{proposition}

Indeed, in the definition of an almost periodic sequence we may
require that each \textsl{prefix} of the sequence has infinitely many
occurences at limited distances (since each substring is a part
of some prefix). And this is guaranteed: any prefix of length
$l<n_s$ appears with period $n_{s+1}$.

The same construction can be explained in a different way.
Consider the positional system where the last digit of integer
$x$ is $x\ \mathrm{mod}\ n_0$, the previous digit is $(x\
\mathrm{div}\ n_0)\ \mathrm{mod}\ n_1$ etc. Then all numbers of
the form $\ldots 0z$ (for any given $z\in [0,n_0)$) are equivalent; we
say that they have rank~$1$. Then we make (for any $y,z$
such that $y\ne 0$) all numbers of the form $\ldots 0 yz$
equivalent and assign rank~$2$ to them, etc.

If the sequence of periods $n_0 < n_1 <n_2<\ldots$ is growing
fast enough, then the equivalence relation does not restrict
significantly the freedom of bit choice: going from left to
right, we see that most of the bits are ``primary'' bits (are
leftmost bits in their equivalence class, not copies of previous
bits; these copies are called ``secondary'' bits, see
Fig.~\ref{bits}).

Indeed, bits of rank $1$ start with $n_0$ primary bits, these
bits are repeated as secondary bits with period $n_1$, so
secondary bits of rank~$1$ form a $n_0/n_1$-fraction of all bits
in the sequence; secondary bits of rank~$2$ form a
$n_1/n_2$-fraction etc. So the sum $\sum_{i}\frac{n_i}{n_{i+1}}$
is the upper bound of the density of ``non-fresh'' bits. More
precise estimate: prefix of any length $N$ has at least $DN$
fresh bits where
š š š š $$
D=\prod_{i} (1-n_i/n_{i+1}).
š š š š $$

This gives a simple proof of Proposition~\ref{prefixes}. For a
given~$\alpha$ choose a computable sequence $n_0<n_1<n_2<\ldots$
that grows fast enough and has $D>\alpha$. Then take a
Martin-L\"of random sequence $\xi$ and place its bits (from left
to right) at all free positions (duplicating bits as required by
the equivalence relation). We get an almost periodic sequence
$\omega$; at least $DN$ bits of $\xi$ can be algorithmically
reconstructed from $\omega$'s prefix of length $N$. It remains
to note that algorithmic transformation cannot increase
complexity and that complexity of $m$-bit prefix of a random
sequence is at least $m-o(m)$ (it would be at least $m$ for
monotone or prefix complexity, but could be $O(\log m)$ smaller
for plain complexity).

\section{Proof of the main result}

Could we apply the same argument (with sequence~$\omega$ from
Proposition~\ref{levin-lemma} instead of a random sequence) to
prove Proposition~\ref{main}? Not directly. To explain the
difficulty and the way to overcome it, consider the simplified
picture where only the equivalence of rank~$1$ is used. Then the
sequence constructed has the form
š š š š $$
\omega=A\,B_0\,A\,B_1\,A\,B_2\,A\,B_3\,A\ldots
š š š š $$
where $A$ is the group of primary bits of rank~$1$ (repeated
with period $n_1$); $A$ and $B_i$ are taken from a sequence
š š š š $$
\xi=A\,B_0\,B_1\,B_2\,B_3\ldots
š š š š $$
(provided by Proposition~\ref{levin-lemma}). If some substring
$x$ of~$\omega$ is located
entirely in $A$ or some $B_i$, its high complexity is
guaranteed by Proposition~\ref{levin-lemma}. However, if $x$
appears on the boundary between $A$ and $B_i$ for some $i>0$,
then $x$ is composed from two substrings of $\xi$ and its
complexity is not guaranteed to be high.

To overcome this difficulty, we need the following stronger
version of Proposition~\ref{levin-lemma}.

\begin{proposition}\label{levin-lemma-modified}
š š š š %
For any $\alpha<1$ there exists a number $N$ and an infinite
binary sequence $\omega$ such that any its substring
š š š š $$
x=\omega_{n}\omega_{n+1}\omega_{n+2}\ldots\omega_{n+k-1}
š š š š $$
of length $k>N$ has high conditional
complexity with respect to previous bits:
š š š š $$
{\mathrm{K}}(\omega_{n}\omega_{n+1}\omega_{n+2}\ldots\omega_{n+k-1}\mid
\omega_{0}\omega_{1}\omega_{2}\ldots\omega_{n-1})
\ge \alpha k.
š š š š $$
\end{proposition}

The proof follows the scheme from~\cite{dls}. Let $\beta<1$ be
greater than~$\alpha$. Let $m$ be some integer number (we will
fix it later). Let the first $m$ bits of~$\omega$ be the
sequence $x$ of length $m$ with maximal prefix complexity
(denoted by ${\mathrm{KP}}$). Then add the next $m$ bits to get the
maximal prefix complexity of the entire sequence. This increase
would be at least $m-O(\log m)$.

[Indeed, for any strings $x$ and $y$ we have
š š š š $$
{\mathrm{KP}}(x,y)={\mathrm{KP}}(x)+
{\mathrm{KP}}(y\mid x,{\mathrm{KP}}(x))+O(1);
š š š š $$
(Kolmogorov -- Levin theorem); if $y$ has been chosen to maximize the
second term in the sum, then ${\mathrm{KP}}(y\mid\ldots)\ge |y|$ and
${\mathrm{KP}}(x,y)\ge{\mathrm{KP}}(x)+|y|-O(1)$. Therefore, for this $y$
š š š š $$
{\mathrm{KP}}(xy)\ge {\mathrm{KP}}(x,y)-{\mathrm{KP}}(|y|)-
O(1)\ge {\mathrm{KP}}(x)+ |y|-O(\log |y|),
š š š š $$
since $(x,y)$ can be reconstructed from $xy$ and $|y|$ and
${\mathrm{KP}}(|y|)=O(\log|y|)$. See~\cite{dls} for details.]

Then we add string $z$ of length $m$ that maximizes ${\mathrm{KP}}(xyz)$
and so on.

In this way we construct a sequence~$\omega=xyz\ldots$ such that
the prefix complexity of its initial segments increases by
$m-c\log m$ for every added block of
$m$ bits. We can choose $m$ such
that $m-c\log m-O(1)>\beta m$.

Then the statement of the Proposition follows from Kolmogorov --
Levin theorem if the substring is ``aligned'' (starts and ends
on the boundaries of length~$m$ blocks). Since $m$ is fixed, the
statement is true for non-aligned blocks of large enough length
(boundary effects are compensated by the difference between $\alpha$
and $\beta$).

Proposition~\ref{levin-lemma-modified} is proven.

Let us explain why this modification helps in the model
situation considered above. If a substring $x$ of the sequence
$AB_0 AB_1 AB_2\ldots$ is on the boundary between $A$ and some
$B_i$, then it can be split into two parts $x_A$ and $x_B$.
The string $x_A$ is a substring of~$A$ and therefore has high
complexity. The string $x_B$ is a substring of some $B_i$ and
therefore also has high complexity and even high conditional
complexity with respect to some prefix containing~$A$.
If we prove that $x_A$ is simple relatively this prifix
we can use Kolmogorov -- Levin theorem to prove that $x$
has high complexity.

Similar arguments work in general case when we have to consider
bits of all ranks. To finish the proof we need the following Lemma:

\textbf{Lemma}. \emph{Let $\omega$ be the sequence satisfying
the statement of Proposition~\ref{levin-lemma-modified}. Then
        $$
\begin{array}{l}
{\mathrm{K}}(V(a_0,b_0),V(a_1,b_1),\ldots,V(a_{s-1},b_{s-1}))\ge\\
\qquad\qquad\alpha L-O(s\log L)-{\mathrm{K}}(a_0\mid a_1)-
{\mathrm{K}}(a_1\mid a_2)-\ldots-{\mathrm{K}}(a_{s-2}\mid a_{s-1})
\end{array}
        $$
for any $a_0<b_0\le a_1<b_1\le\dots\le a_{s-1}<b_{s-1}$, where
$V(a,b)$ stands for $\omega_a \omega_{ a+1}\dots \omega_{b-1}$
and $L=(b_0-a_0)+(b_1-a_1)+\dots+(b_{s-1}-a_{s-1})$.}

In fact, for Proposition~\ref{main} we need only the case $s=3$
of this Lemma.

The proof of Lemma is based on Kolmogorov -- Levin theorem about
complexity of pairs. The statement of Proposition~\ref{levin-lemma-modified}
guarantees the following inequality:
        $$
{\mathrm{K}}(V(a_{s-1},b_{s-1})\mid V(0,a_{s-1}))\ge
\alpha(b_{s-1}-a_{s-1})-O(\log L).\eqno{(*)}
        $$
We will prove the following inequality of any $i=0,1,\dots,s-2$:
        $$
\begin{array}{r}
{\mathrm{K}}(V(a_i,b_i),V(a_{i+1},b_{i+1}),\ldots,V(a_{s-1},b_{s-1})\mid V(0,a_i))-\qquad\qquad\\
{\mathrm{K}}(V(a_{i+1},b_{i+1}),\ldots,V(a_{s-1},b_{s-1})\mid V(0,a_{i+1}))\ge\qquad\\
\alpha(b_i-a_i)-O(\log L)-{\mathrm{K}}(a_i\mid a_{i+1}).
\end{array}\eqno{(**)}
        $$
If we add up (**) for all $i=0,1,\dots,s-2$ with (*) we obtain
the required inequality (and even stronger one
with relative complexity in the left-hand side).
Let us prove the inequality~(**) now. By $W$ we denote the
sequence $(V(a_{i+1},b_{i+1}),\ldots,V(a_{s-1},b_{s-1}))$. The
following inequality follows from the Kolmogorov -- Levin theorem
and the statement of Proposition~\ref{levin-lemma-modified}:
        $$
\begin{array}{r}
{\mathrm{K}}(V(a_i,b_i),W\mid V(0,a_i))-
{\mathrm{K}}(W\mid V(0,a_i),V(a_i,b_i))=\qquad\qquad\\
{\mathrm{K}}(V(a_i,b_i)\mid V(0,a_i))-O(\log L)\ge\alpha(b_i-a_i)-O(\log L).
\end{array}
        $$
To finish the proof of Lemma, let us prove the inequality
        $$
{\mathrm{K}}(W\mid V(0,a_{i+1}))\le{\mathrm{K}}
(W\mid V(0,a_i),V(a_i,b_i))+{\mathrm{K}}(a_i\mid a_{i+1})+O(\log L).
        $$
One can obtain~$W$ from $V(0,a_{i+1})$ in the following way:
find~$a_{i+1}$ using the length of the string $V(0,a_{i+1})$,
convert~$a_{i+1}$ into~$a_i$ by the shortest program,
compute~$b_i$ by adding difference~$b_i-a_i$ to~$a_i$,
cut intervals $[0,a_i)$ and $[a_i,b_i)$ from string $V(0,a_{i+1})$
and execute the shortest program that converts
$(V(0,a_i),V(a_i,b_i))$ into~$W$.
This needs ${\mathrm{K}}(W\mid V(0,a_i),V(a_i,b_i))+
{\mathrm{K}}(a_i\mid a_{i+1})+O(\log L)$~bits
to obtain~$W$ from $V(0,a_{i+1})$. The inequality is proven, q.e.d.

The proof of Proposition~\ref{main} uses the same construction
as proof of Proposition~\ref{prefixes} but it takes a sequence
satisfying the statement of Proposition~\ref{levin-lemma-modified}
instead of a random sequence.

Let $v$~be a sequence satisfying the statement
of Proposition~\ref{levin-lemma-modified} with some $\alpha'>\alpha$ and $\omega$~be
the resulting sequence (if we apply the construction
of an almost periodic sequence to the sequence~$v$).
It has been proved before that $\omega$~is an almost
periodic sequence. We need only to prove the following
estimate of a complexity of any substring of~$\omega$:
        $$
{\mathrm{K}}(\omega_{m}\omega_{m+1}\omega_{m+2}\ldots\omega_{m+k-1})\ge\alpha k.
        $$
for any sufficiently long~$k$ and for any~$m$.

Suppose that sequence~$\{n_j\}$ grows fast enough, i.e.
$\sum_{j=1}^{\infty}\frac{n_{j-1}}{n_j}<\frac{\alpha'-\alpha}{2}$.
Suppose $i$~is the smallest index such that~$n_i\ge k$.
Due to our construction of sequence~$\omega$ any element
of~$\omega$ corresponds to some element of~$v$. Different elements
of $\omega_{m}\omega_{m+1}\ldots\omega_{m+k-1}$ of rank
not less than~$i$ (i.e. elements repeated with period~$n_i$ or greater
by our construction) correspond to different elements of~$v$ because
the distance between elements of the given substring of~$\omega$ is
less than~$n_i$ (and less than the period). It is easy to prove that
in this substring the density of elements of small rank (less than~$i$)
is not greater than $\alpha'-\alpha$.

Indeed, the number of elements of rank~$j$ on any interval of length~$n_j$
is equal to $n_{j-1}$ and we can cover the given interval
of length~$k$ with at most $\frac{k}{n_j}+1$ intervals of length~$n_j$.
Therefore the number of elements of rank~$j$ on the given interval
is not greater than $n_{j-1}\left(\frac{k}{n_j}+1\right)$.
So the density of elements of rank less than~$i$
in the given substring is not greater than
$\sum_{j=1}^{i-1}\left(\frac{n_{j-1}}{n_j}+\frac{n_{j-1}}{k}\right)\le
2\sum_{j=1}^{i-1}\frac{n_{j-1}}{n_j}<\alpha'-\alpha$ due to
our assumption about growing of~$\{n_j\}$.

Hence the substring $\omega_{m}\ldots\omega_{m+k-1}$ corresponds
to some intervals in~$v$. Throw away all elements of small ranks
from these intervals of~$v$ and denote the remaining intervals
by $[a_0,b_0),\dots,[a_{s-1},b_{s-1})$, where $a_0<b_0\le\dots\le a_{s-1}<b_{s-1}$.
The number of these intervals is at most~$3$. Indeed, we can enumerate all
elements of $\omega_{m}\ldots\omega_{m+k-1}$  from left to right, not counting
elements of small ranks, and for each element find the corresponding element of~$v$.
The index of corresponding element will increase by~$1$ every time except when we
cross a point of type~$n_i j$ or $n_i j+n_{i-1}$ (where $j$ is integer).
But there are at most $2$~points of this type in the interval of length~$k$
so there are at most $3$~corresponding intervals.

Substrings~$V(a_0,b_0),\dots,V(a_{s-1},b_{s-1})$ (defined as in Lemma)
can be computed by an algorithm using the given substring of~$\omega$.
The algorithm needs only to know the value of $m\mathrel{\mathrm{mod}}n_{i-1}$
for finding elements with small rank (less than~$i$) and the relative positions
of elements of~$\omega_{m}\ldots\omega_{m+k-1}$ corresponding to $v_{a_j}$
and $v_{b_j-1}$ where $j=0,1,\dots,s-1$. Because $s\le 3$ only
a logarithmical amount of additional bits is needed. So we can prove the following
inequality to finish the proof of Proposition~\ref{main}:
        $$
{\mathrm{K}}(V(a_0,b_0),\dots,V(a_{s-1},b_{s-1}))\ge\alpha k-O(\log k).
        $$
We can use Lemma for this because $\alpha' L>\alpha k$,
where $L=(b_0-a_0)+(b_1-a_1)+\dots+(b_{s-1}-a_{s-1})$
(we have already proved that in this substring the density of
elements of small rank is not greater than $\alpha'-\alpha$,
hence $k-L\le(\alpha'-\alpha)k$).

If we prove that ${\mathrm{K}}(a_j\mid a_{j+1})=O(\log k)$
we will finish the proof of the proposition.
Suppose we know~$a_{j+1}$. We can find~$a_j$ in the following way.
Find the element of the given substring of~$\omega$ corresponding
to~$v_{a_{j+1}}$. Add to the index of the found element the difference between
the indexes of the elements of the given substring corresponding to~$v_{a_j}$
and~$v_{a_{j+1}}$ (this difference is not greater than the lenght of the
given substring, i.e., we use only a logarithmical amount of memory).
We get an element of~$\omega$ corresponding to~$v_{a_j}$. It can be used
to calculate~$a_j$. But the first step of this algorithm uses knowing
the position of the given substring which needs an unlimited amount of memory.
We can avoid using this position if we notice that the rank~$i$ of
elements of~$\omega$ corresponding to~$v_{a_j}$ is not greater than the rank~$I$ of
elements of~$\omega$ corresponding to~$v_{a_{j+1}}$ (because $a_j<a_{j+1}$).
So $n_I$~is a multiple of~$n_i$. Hence at the first step we can take any
element of~$\omega$ corrensponding to~$v_{a_{j+1}}$ (for example, the first one).
We get the same result since the elements corresponding
to~$v_{a_j}$ repeat with period~$n_i$ and the elements corresponding
to~$v_{a_{j+1}}$ repeat with period~$n_I$.

Therefore we construct the algorithm proving
that~${\mathrm{K}}(a_j\mid a_{j+1})=O(\log k)$, and so the proof of
Proposition~\ref{main} is complete.

\textbf{Remarks}.

1. Proposition~\ref{main} implies the existence of a
bi-infinite almost periodic sequence with complex substrings
(using the standard compactness argument; this argument
can be even
simplified for the special case of almost periodic sequences).

2. The proof of Proposition~\ref{main} works for relativized
version of complexity. Therefore we get (as explained above) the
following (pure combinatorial) strong version of
Proposition~\ref{forbidden-strings}:

\textbf{Corollary}. \emph{Assume that for each $n$ a set~$F_n$
of forbidden substrings of length~$n$ is fixed, and the size
of~$F_n$ is at most $2^{\alpha n}$.
Then there exists an infinite almost periodic binary sequence
$\omega$ and a constant $N$ such that for any $n>N$ the sequence
$\omega$ does not have a substring $x$ that belongs to $F_n$.}

\section{Multidimensional case}\label{multidimensional}

Similar but more delicate arguments could be applied to
multidimensional case too.

A $d$-dimensional sequence $\omega: \bbbz^d\to\{0,1\}$ is
\emph{almost periodic} if for any cube~$x$ that appears
in~$\omega$ there exists a number~$k$ such that any subcube with
side~$k$ contains $x$ inside.

\begin{proposition}\label{multidimensional-main}
š š š š %
Fix an integer $d\ge 1$. Let $\alpha$ be a positive number less
than~$1$. There exists an almost periodic $d$-dimensional
sequence~$\omega$ such that any sufficiently large subcube~$x$
of~$\omega$ has large complexity:
š š š š $$
{\mathrm{K}}(x)\ge\alpha\cdot\mbox{volume}(x)
š š š š $$
\end{proposition}

Here volume is the number of points, i.e., $\mbox{side}^d$.

In the multidimensional case the complexity argument needs
Proposition~\ref{levin-lemma-modified} even if we do not insist
that $\omega$ is almost periodic.

Informally, the idea of the proof can be explained as follows.
Consider, for example, the case $d=2$. Take a sequence $v$
from Proposition~\ref{levin-lemma-modified} and write down its
terms along a spiral.

\begin{center}
\includegraphics[scale=1]{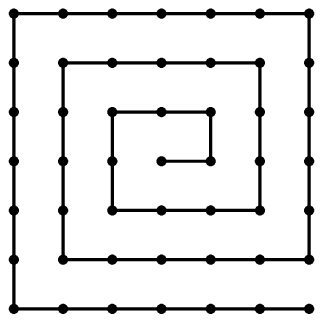}
\end{center}

Then we need to bound the complexity of a cube (i.e., square).
This square contains several substrings of the sequence~$v$.
(Unlike the previous case where only $3$ substrings were needed,
now the number of substrings is proportional to the side of the
square.) Then we apply the Lemma to these substrings to get the
bound for the complexity of the entire square.

This works if we do not require $\omega$ to be almost periodic
(so the argument above could replace the combinatorial proof
using Lovasz lemma). It needs additional modifications to get
the almost periodic sequence. Similar to one-dimensional
construction, the cube $[-n_0,n_0)^d$ is duplicated periodically
in all directions with shifts being multiples of $n_1$
(where $n_0\mid n_1$); the cube $[-n_1,n_1)^d$ is duplicated
with shifts being multiples of $n_2$ (where $n_1\mid n_2$), etc.

As in one-dimensional case, it is easy to see that this
construction guarantees that $\omega$ is almost periodic.
Let $v$~be a sequence satisfying the statement
of Proposition~\ref{levin-lemma-modified} with
some $\alpha'>\alpha$. We sort all new positions
of~$\omega$ by rank (the element has rank~$j$
if it is duplicated with period~$n_j$ by the structure
described) then by coordinated in lexicographical order.
Then we fill the positions with the elements of~$v$
in this order. Let $B=[m_1,m_1+k)\times[m_2,m_2+k)\times\dots
\times[m_d,m_d+k)$ is a cube. We need only to prove that
cube~$B$ in the sequence~$\omega$ has high comlexity:
        $$
{\mathrm{K}}(\omega_B)\ge\alpha k^d.
        $$
\begin{figure}[h]
\begin{center}
\includegraphics[scale=0.75]{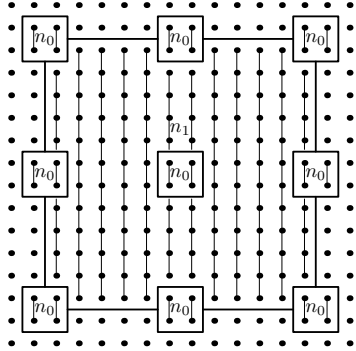}
\end{center}
\caption{Duplicated cubes in two dimentional case.}
\end{figure}

Suppose that sequence~$\{n_j\}$ grows fast enough, i.e.
$\sum_{j=1}^{\infty}\frac{n_{j-1}}{n_j}<\frac{\alpha'-\alpha}{4}$.
Suppose $i$~is the smallest index such that~$n_i\ge k$.
Due to our construction of sequence~$\omega$ any element
of~$\omega$ corresponds to some element of~$v$. Different elements
of~$\omega$ of rank not less than~$i$ in cube~$B$
correspond to different elements of~$v$ because
the distance between elements of the given cube is
less than~$n_i$ (and less than the period). It is easy to prove that
in this cube the density of elements of small rank (less than~$i$)
is not greater than $\alpha'-\alpha$.

Indeed, the number of elements of rank~$j$ on any vertical
(i.e., parallel to the last axis) interval of length~$n_j$
is zero or $2n_{j-1}$ and we can cover the given cube
of side~$k$ with at most $k^{d-1}(\frac{k}{n_j}+1)$ vertical
intervals of length~$n_j$. Therefore the number
of elements of rank~$j$ on cube~$B$
is not greater than $2n_{j-1}k^{d-1}\left(\frac{k}{n_j}+1\right)$.
So the density of elements of rank less than~$i$
in the given cube is not greater than
$2\sum_{j=1}^{i-1}\left(\frac{n_{j-1}}{n_j}+\frac{n_{j-1}}{k}\right)\le
4\sum_{j=1}^{i-1}\frac{n_{j-1}}{n_j}<\alpha'-\alpha$ due to
our assumption about growing of~$\{n_j\}$.

Hence cube~$B$ corresponds to some intervals in~$v$.
Throw away all elements of small ranks
from these intervals of~$v$ and denote the remaining intervals
by $[a_0,b_0),\dots,[a_{s-1},b_{s-1})$,
where $a_0<b_0\le\dots\le a_{s-1}<b_{s-1}$.
The number of these intervals is at most~$4k^{d-1}$.
Indeed, we can enumerate all elements of each vertical
interval of length~$k$ in our cube from bottom to top
(from small last coordinate to big one), not counting
elements of small ranks, and for each element find the corresponding
element of~$v$. The index of corresponding element will increase
by~$1$ every time except when we cross a point of
type~$n_i j$, $n_i j+n_{i-1}$ or $n_i j-n_{i-1}$ (where $j$ is integer).
But there are at most $3$~points of this type in any vertical interval
of length~$k$ so there are at most $4$~corresponding intervals
for each vertical interval. But the number of vertical intervals
of length~$k$ in cube~$B$ is equal to $k^{d-1}$, so the total
number of corresponding intervals $s\le 4k^{d-1}$.

Substrings~$V(a_0,b_0),\dots,V(a_{s-1},b_{s-1})$ (defined as in Lemma)
can be computed by an algorithm using the given substring of~$\omega$.
The algorithm needs only to know the value of
$m_j\mathrel{\mathrm{mod}}n_{i-1}$, where $j=1,2,\dots,d$,
for finding elements with small rank (less than~$i$) and
the relative positions in the cube~$B$ corresponding to $v_{a_j}$
and $v_{b_j-1}$ where $j=0,1,\dots,s-1$. Because $s\le 4k^{d-1}$
the algorithm needs only $O(k^{d-1}\log k)$ bits.
So we can prove the following inequality to finish
the proof of Proposition~\ref{multidimensional-main}:
        $$
{\mathrm{K}}(V(a_0,b_0),\dots,V(a_{s-1},b_{s-1}))\ge
\alpha k^d-O(k^{d-1}\log k)
        $$
(the value~$O(k^{d-1}\log k)$ is compensated by a small change of~$\alpha$).
We can use Lemma for this because $\alpha' L>\alpha k$,
where $L=(b_0-a_0)+(b_1-a_1)+\dots+(b_{s-1}-a_{s-1})$
(we have already proved that in this cube the density of
elements of small rank is not greater than $\alpha'-\alpha$,
hence $k-L\le(\alpha'-\alpha)k$).

If we prove that ${\mathrm{K}}(a_j\mid a_{j+1})=O(\log k)$
we will finish the proof of the proposition.
Suppose we know~$a_{j+1}$. We can find~$a_j$ in the following way.
Find some element of~$\omega$ corresponding to~$v_{a_{j+1}}$
(for example, the smallest one).
Add to the index of the found element the difference between
the positions in the given cube corresponding to~$v_{a_j}$
and~$v_{a_{j+1}}$ (this difference is not greater than the side of the
cube, i.e., we use only a logarithmical amount of memory).
We get an element of~$\omega$ corresponding to~$v_{a_j}$.
It can be used to calculate~$a_j$. This can be proven
the same way as in Proposition~\ref{main}. If at the first step
we found the element in cube~$B$ corresponding to~$v_{a_{j+1}}$
we obviously would get~$v_{a_j}$ as a result. Notice that
the rank~$i$ of elements of~$\omega$ corresponding to~$v_{a_j}$
is not greater than the rank~$I$ of elements of~$\omega$ corresponding
to~$v_{a_{j+1}}$ (because $a_j<a_{j+1}$). So $n_I$~is a multiple of~$n_i$
and the result does not depend on the element corresponding to~$v_{a_{j+1}}$
since the elements corresponding to~$v_{a_j}$ repeat with period~$n_i$
and the elements corresponding to~$v_{a_{j+1}}$ repeat with period~$n_I$.

Therefore we construct the algorithm proving
that~${\mathrm{K}}(a_j\mid a_{j+1})=O(\log k)$, and so the proof of
Proposition~\ref{multidimensional-main} is complete.

\section{Remarks}

Kolmogorov complexity is often used in combinatorial
constructions as the replacement of counting arguments. (Instead
of proving that the total number of objects is larger that the
number of ``bad'' objects we prove that an object of maximal
complexity is ``good''.) Sometimes people even say that the use
of Kolmogorov complexity is just a simple reformulation that often
hides the combinatorial essence of the argument.

In our opinion this is not always true. Even without the almost
periodicity requirement the two natural proofs of
Proposition~\ref{forbidden-strings} (using complexity argument
and Lovasz lemma) are quite different. The proof
of Proposition~\ref{levin-lemma} uses prefix complexity and cannot be
directly translated into a counting argument. On the other hand,
the use of Lovasz lemma in a combinatorial proof cannot be
easily reformulated in terms of Kolmogorov complexity.
(Moreover, for almost periodic case we don't know how to apply
Lovasz lemma argument and complexity proof remains the only one
known to us.)

\section{Acknowledgements}

The authors would like to thank Alexander Shen and Yury Pritykin for help.

\end{document}